\documentclass[journal]{IEEEtran}
\usepackage{cite}
\usepackage[pdftex]{graphicx}
\usepackage{enumitem}
\usepackage{amsfonts}
\usepackage{amsmath}
\usepackage{color}
\usepackage{multirow}
\newtheorem{thm}{Theorem}[section]
\newtheorem{lemma}{Lemma}[section]
\newtheorem{defin}{Definition}[section]
\newtheorem{assume}{Assumption}[section]
\usepackage[ruled,norelsize]{algorithm2e}

\makeatletter
\newcommand{\removelatexerror}{\let\@latex@error\@gobble}
\makeatother

\begin{document}

\title{Distributionally Robust Chance Constrained Optimal Power Flow Assuming Unimodal Distributions with Misspecified Modes}

\author{Bowen~Li,~\IEEEmembership{Student Member,~IEEE,}
        Ruiwei~Jiang,~\IEEEmembership{Member,~IEEE,}
        and~Johanna~L.~Mathieu,~\IEEEmembership{Member,~IEEE}% <-this % stops a space
\thanks{This work is supported by the U.S. National Science Foundation Awards CCF-1442495 and CMMI-1662774. B. Li and J. L. Mathieu are with the Department of Electrical Engineering and Computer Science, University of Michigan at Ann Arbor, Ann Arbor, MI 48109 USA (e-mail: libowen@umich.edu; jlmath@umich.edu). R. Jiang is with the Department of Industrial and Operations Engineering, University of Michigan at Ann Arbor, Ann Arbor, MI 48109 USA (e-mail: ruiwei@umich.edu).}% <-this % stops a space
}

% \markboth{Submitted to IEEE transactions on Control of Network systems}
% {Li \MakeLowercase{\textit{et al.}}: Distributionally Robust Optimal Power Flow Assuming Ambiguous-Mode Unimodality}

% make the title area
\maketitle

% As a general rule, do not put math, special symbols or citations
% in the abstract or keywords.
\begin{abstract}
Chance constrained optimal power flow (CC-OPF) formulations have been proposed to minimize  operational costs while controlling the risk arising from uncertainties like renewable generation and load consumption. To solve CC-OPF, we often need access to the (true) joint probability distribution of all uncertainties, which is rarely known in practice. A solution based on a biased estimate of the distribution can result in poor reliability. To overcome this challenge, recent work has explored distributionally robust chance constraints, in which the chance constraints are satisfied over a family of distributions called the ambiguity set. Commonly, ambiguity sets are only based on moment information (e.g., mean and covariance) of the random variables;
however, specifying additional characteristics of the random variables reduces conservatism and cost. Here, we consider ambiguity sets that additionally incorporate unimodality information. In practice, it is difficult to estimate the mode location from the data and so we allow it to be potentially misspecified. We formulate the problem and derive a separation-based algorithm to efficiently solve it. Finally, we evaluate the performance of the proposed approach on a modified IEEE-30 bus network with wind uncertainty and compare with other distributionally robust approaches. We find that a misspecified mode significantly affects the reliability of the solution and the proposed model demonstrates a good trade-off between cost and reliability.
\end{abstract}

% Note that keywords are not normally used for peerreview papers.

% For peer review papers, you can put extra information on the cover
% page as needed:
% \ifCLASSOPTIONpeerreview
% \begin{center} \bfseries EDICS Category: 3-BBND \end{center}
% \fi
%
% For peerreview papers, this IEEEtran command inserts a page break and
% creates the second title. It will be ignored for other modes.
\IEEEpeerreviewmaketitle

\section{Introduction}
% The very first letter is a 2 line initial drop letter followed
% by the rest of the first word in caps.
%
% form to use if the first word consists of a single letter:
% \IEEEPARstart{A}{demo} file is ....
%
% form to use if you need the single drop letter followed by
% normal text (unknown if ever used by the IEEE):
% \IEEEPARstart{A}{}demo file is ....
%
% Some journals put the first two words in caps:
% \IEEEPARstart{T}{his demo} file is ....
%
% Here we have the typical use of a "T" for an initial drop letter
% and "HIS" in caps to complete the first word.
\IEEEPARstart{W}{ith} higher penetrations of renewable generation, uncertainties have increasing influence on power system operation and hence need to be carefully considered in scheduling problems, such as  optimal power flow (OPF). To manage the risk arising from uncertainties, different stochastic OPF approaches have been studied. Among these formulations, CC-OPF has been proposed to directly control the constraint violation probability below a pre-defined threshold \cite{cc1,cc2,vrakopoulou_probabilistic_2013,bienstock_chance_2014,cc4,mariatsg,bowentsg}. Traditional methods to solve chance constrained programs require knowledge of the joint probability distribution of all uncertainties, which may be unavailable or inaccurate. However, biased estimate may yield poor out-of-sample performance. Randomized techniques such as scenario approximation \cite{campi2009, Margellos2014}, which provides a priori guarantees on reliability, require the constraints to be satisfied over a large number of uncertainty samples. The solutions from these approaches are usually overly conservative with high costs \cite{DRyiling,bowentsg}. Another popular approach is to assume that the uncertainties follow a parametric distribution such as Gaussian \cite{bienstock_chance_2014,cc4,bowentsg}. The resulting CC-OPF is often easier to solve but the solution may have low reliability unless the assumed probability distribution happens to be close to the true one.

As an alternative, distributionally robust chance constrained (DRCC) OPF models do not depend on a single estimate of the probability distribution \cite{dr3,dr4,summers,dr5,bowenpscc, unimodalBL,roald1,DRyiling}. More specifically, DRCC models consider a family of distributions, called the ambiguity set, that share certain statistical characteristics and requires that the chance constraint holds with respect to all distributions within the ambiguity set \cite{el2003worst,delage2010distributionally,stellato2014data,ruiwei}. %These statistical properties are generally computed from empirical data and hence affected by sampling errors. In general, the estimated results from different data samples will distribute around a certain region. If we use any single value within the region to represent the true property, the resulting probability distribution might be biased and result in low-reliability solutions. Specifically, for unimodality, ambiguity on mode location affects the skewness and the tail tendency of the distribution. Accurate mode estimation will help reduce the objective value in a DR optimization problem and the resulted optimal solutions will satisfy our probability requirement on constraint violation. Hence, to avoid the inaccurate mode estimate, we would like to consider ambiguous modes in the ambiguity sets.
Most existing work characterizes the ambiguity set based on moment information obtained from historical data of the uncertainty (see, e.g.,~\cite{DRyiling,dr3,dr5,summers}). For example, a commonly adopted ambiguity set consists of all distributions whose mean and covariance agree with their corresponding sample estimates~\cite{DRyiling,dr3,summers}. Many uncertainty distributions (e.g., those associated with wind forecast error) are unimodal and so, recently, unimodality has been incorporated to strengthen the ambiguity set and reduce the conservatism of DRCC models \cite{summers,unimodalBL,roald1}. However, as compared to the moments, the mode location is more likely to be misspecified in sample-based estimation.

In this paper, we study a DRCC model with an ambiguity set based on moment and unimodality information with a potentially misspecified mode location. To the best of our knowledge, this paper is the first work discussing misspecification of a value related to a structural property, though others have considered misspecification of moments \cite{delage2010distributionally,dr5,el2003worst,ruiwei,DRyiling} and misspecification of distributions \cite{dr4,DRbitar}.
Our main theoretical result shows that the distributionally robust chance constraints can be recast as a set of second-order conic (SOC) constraints. Furthermore, we derive an iterative algorithm to accelerate solving the reformulation. In this algorithm, we begin with a relaxed formulation, and in each iteration, we efficiently find the most violated SOC constraint, if any, or terminate with a globally optimal solution. We apply the theoretical results to a direct current (DC) OPF problem and conduct a case study using a modified IEEE 30-bus system with wind power. We compare our results (operational cost, reliability, computational time, and optimal solutions) to those obtained using four alternative ambiguity sets \cite{unimodalBL,stellato2014data,roald1,DRyiling}.

The remainder of this paper is organized as follows. Section~\ref{sec: winderror} empirically verifies the (multivariate) unimodality of wind forecast errors and explores misspecification of the mode location. The proposed DRCC model and ambiguity set are introduced in Section~\ref{sec: drcc} and the main theoretical results are presented in Section~\ref{sec:main_drcc}. Section~\ref{sec:casestudy} includes the case studies and Section~\ref{sec:summary} concludes the paper.

\section{Unimodality of Wind Forecast Errors \& Error in Mean and Mode Estimates}\label{sec: winderror}
\begin{figure}
\centering
\vspace{.3cm}
\includegraphics[width=3.5in]{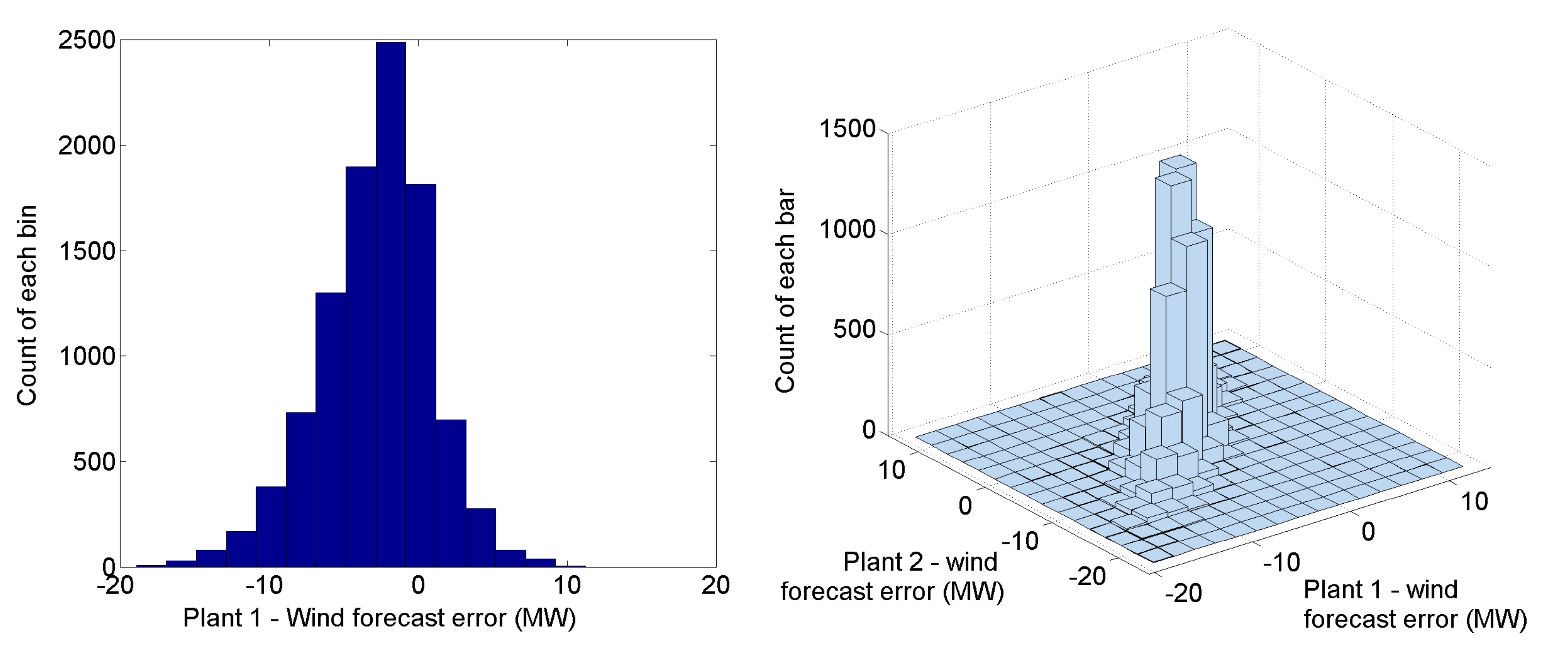}
\caption{Histograms of univariate and bivariate wind forecast errors ($15$ bins).}
\vspace{-.3cm}
\label{fig:windunimodal}
\end{figure}

In this section, we first empirically verify the unimodality of wind forecast error distributions using 10,000 data samples from \cite{mariatsg,bowentsg} with statistical outliers omitted (total probability $<0.1\%$). The samples were generated using a Markov Chain Monte Carlo mechanism \cite{MCMC} based on real data that includes both hourly forecast and actual wind generation in Germany. In Fig.~\ref{fig:windunimodal}, we depict the histograms of univariate and bivariate wind forecast errors with $15$ bins. Both histograms empirically justify our assumption that the probability distribution of wind forecast errors is unimodal.

\begin{figure}
\centering
\vspace{.3cm}
\includegraphics[width=3.5in]{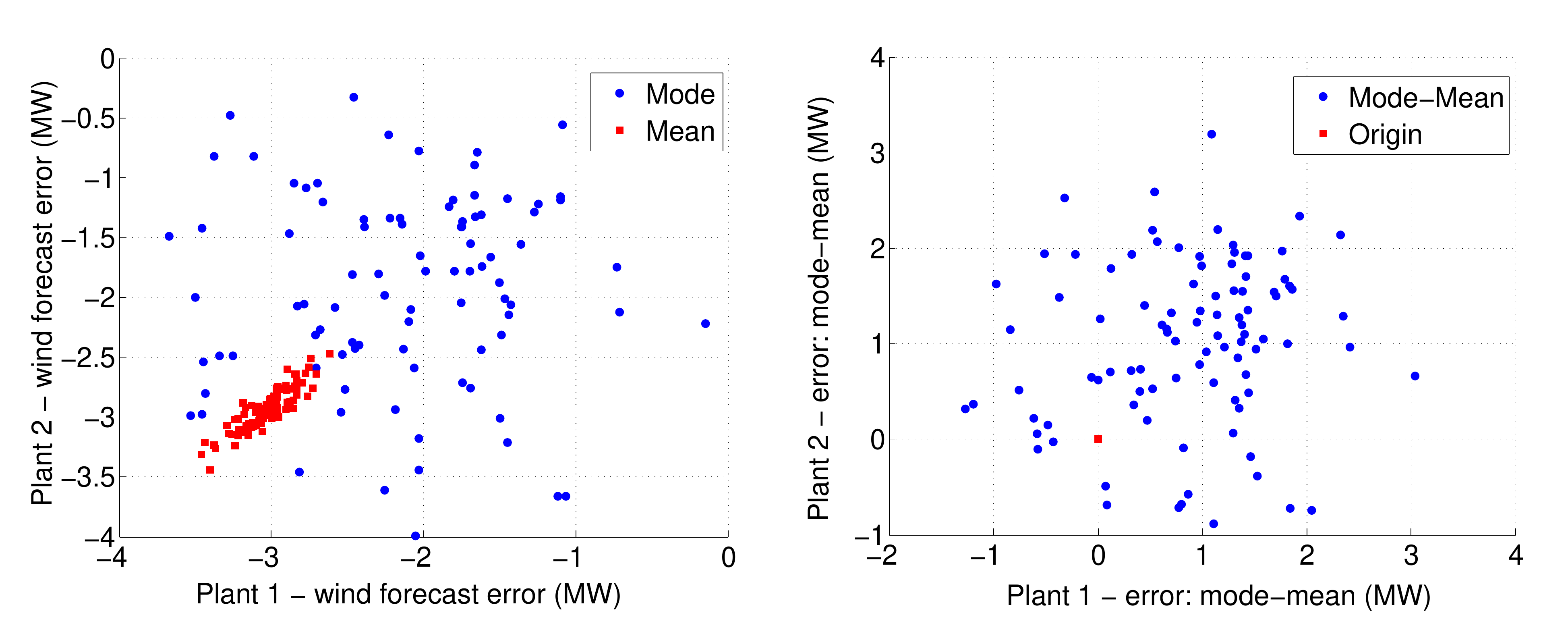}
\caption{Scatter plots of mode and mean estimates from data samples (left) and mode vs. mean differences (right).}
\vspace{-.3cm}
\label{fig:windscatter}
\end{figure}

Next, we empirically evaluate the errors of mean and mode estimates (i.e., the peak location in the histogram). We randomly extract 100 groups of samples, each group containing 500 data points, from the wind forecast error data pool. For each group of samples, we estimate the mean by taking sample averages and estimate the mode by identifying the center of the highest bin in the 15-bin histogram. In Fig.~\ref{fig:windscatter}, we plot all the mean and mode estimates and the differences between them. From the left subfigure, we observe that sampling errors have larger impacts on mode estimates than on mean estimates. From the right subfigure, we observe that the mode estimate can deviate from the corresponding mean estimate in all directions. This indicates the importance of considering the misspecification of mode location in DRCC models, because the mode-mean deviation shows the skewness of the uncertainty. As a result, if we misspecify the mode location (e.g., by modeling a right-skewed distribution as a left-skewed one, see Section~\ref{sec: numexp} for an example), then we may mistakenly relax the chance constraint and get poor out-of-sample performance.

\section{DRCC Formulation}\label{sec: drcc}
\subsection{General Formulation}
In this paper, we consider the following physical constraint under uncertainty:
\begin{equation}
a(x)^{\top}\xi\leq b(x), \label{eq: cons}
\end{equation}
where $x\in\mathbb{R}^l$ represents an $l$-dimensional decision variable, and $a(x):\mathbb{R}^l\to\mathbb{R}^n$ and $b(x):\mathbb{R}^l\to\mathbb{R}$ represent two affine functions of $x$. Uncertainty $\xi\in\mathbb{R}^n$ represents an $n$-dimensional random vector defined on probability space $(\mathbb{R}^n,\mathcal{B}^n,\mathbb{P}_{\xi})$ with Borel $\sigma$-algebra $\mathcal{B}^n$ and probability distribution $\mathbb{P}_{\xi}$. The assumption that $a(x)$ and $b(x)$ are affine in $x$ is a standard assumption in existing DRCC models and consistent with the DRCC DC OPF.

To manage constraint violations due to uncertainty, one natural way is to ensure that \eqref{eq: cons} is satisfied with at least a pre-defined probability threshold $1-\epsilon$, which leads to the following chance constraint \cite{charnes1958cost, miller1965chance}:
\begin{equation}
\mathbb{P}_{\xi}\left(a(x)^{\top}\xi\leq b(x)\right)\geq1-\epsilon,\label{eq: cc}
\end{equation}
where $1-\epsilon$ normally takes a large value (e.g., $0.99$).

\subsection{Distributionally Robust Formulation}
In reality, it may be challenging to access the (true) joint probability distribution $\mathbb{P}_{\xi}$. Oftentimes we may only have a set of historical data and certain domain knowledge of $\xi$. In this case, we can consider the following distributionally robust chance constraint:
\begin{equation}
\inf_{\mathbb{P}_{\xi}\in \mathcal{D}_\xi}\mathbb{P}_{\xi}\left(a(x)^{\top}\xi\leq b(x)\right)\geq1-\epsilon.\label{eq: drcc}
\end{equation}
Instead of assuming that $\mathbb{P}_{\xi}$ takes a specific form, we consider an ambiguity set $\mathcal{D}_\xi$ consisting of plausible candidates of $\mathbb{P}_{\xi}$. Then, we require that chance constraint \eqref{eq: cc} holds with respect to all distributions in $\mathcal{D}_{\xi}$.

\subsection{Ambiguity Sets}
In this paper, we consider three ambiguity sets, denoted as $\mathcal{D}_{\xi}^i$ for $i=1,2,3$, that are defined by a combination of moment and unimodality information. Precisely, we consider a generalized notion of unimodality defined as follows.

\vspace{.2cm}
\begin{defin}{($\alpha$-Unimodality~\cite{unimodal})}\label{def: unimodal}
For any fixed $\alpha > 0$, a probability distribution $\mathbb{P}$ on $\mathbb{R}^n$ is called $\alpha$-unimodal with mode $0$ if $t^\alpha\mathbb{P}(B/t)$ is non-decreasing in $t>0$ for every Borel set $B\in\mathcal{B}^n$.
\end{defin}
\vspace{.2cm}

From the definition, we notice that $\alpha$ parameterizes the ``degree of unimodality.'' When $\alpha=n=1$, the definition coincides with the classical univariate unimodality with mode $0$. When $\alpha=n>1$, the density function of $\xi$ (if exists) peaks at the mode and is non-increasing in any directions moving away from the mode. As $\alpha\to\infty$, the requirement of unimodality gradually relaxes and eventually vanishes. Under Definition \ref{def: unimodal}, we define the following three ambiguity sets:
\vspace{.1cm}
 {\bf Ambiguity set 1:} (moment information only)
\begin{align}
\mathcal{D}_{\xi}^1 := \left\{\mathbb{P}_{\xi} \in \mathcal{P}^n: \mathbb{E}_{\mathbb{P}_{\xi}}[\xi] = \mu, \ \mathbb{E}_{\mathbb{P}_{\xi}}[\xi \xi^{\top}] = \Sigma\right\}, \label{Dset1}
\end{align}
 {\bf Ambiguity set 2:} (moment and $\alpha$-unimodality, fixed mode)
\begin{align}
&\mathcal{D}_{\xi}^2 := \bigl\{\mathbb{P}_{\xi} \in \mathcal{P}^n_{\alpha} \cap \mathcal{D}_{\xi}^1: \ \mathcal{M}(\xi)=m_t \bigr\}, \label{Dset2}
\end{align}
{\bf Ambiguity set 3:} (moment and $\alpha$-unimodality, misspecified mode)
\vspace{-.3cm}
\begin{align}
&\mathcal{D}_{\xi}^3 := \bigl\{\mathbb{P}_{\xi} \in \mathcal{P}^n_{\alpha} \cap \mathcal{D}_{\xi}^1: \ \mathcal{M}(\xi)\in\Xi\bigr\}, \label{Dset3}
\end{align}
where $\mathcal{P}^n_{\alpha}$ and $\mathcal{P}^n$ denote all probability distributions on $\mathbb{R}^n$ with and without the requirement of $\alpha$-unimodality respectively; $\mu$ and $\Sigma$ denote the first and second moments of $\xi$; and $\mathcal{M}(\xi)$ denotes a function returning the true mode location of $\xi$ with $m_{t}$ and $\Xi$ representing a single mode value and a connected and compact set. The compact set can be constructed using possible mode estimates calculated from samples of historical data.

Among these three ambiguity sets, we use $\mathcal{D}_{\xi}^1$ as a benchmark.  Set $\mathcal{D}_{\xi}^2$ is a special case of $\mathcal{D}_{\xi}^3$, i.e.,  $\Xi$ only contains a single value $m_t$. In practice, since the mode estimate is influenced by sampling errors, the mode estimates from data samples are not the same single values but distribute around a certain area. The shape of this area decides the underlying structural skewness in the uncertainty distribution. Hence, we compare $\mathcal{D}_{\xi}^2$ and $\mathcal{D}_{\xi}^3$ to see how misspecified mode estimates affect the DRCC problem. In this paper, we do not additionally consider misspecified moments since this topic has been well-studied \cite{delage2010distributionally,dr5,el2003worst,ruiwei} and our main results can be easily extended based on these existing works.

% needed in second column of first page if using \IEEEpubid
%\IEEEpubidadjcol
\subsection{Numerical Example}\label{sec: numexp}
\begin{figure}
\centering
\vspace{.3cm}
\includegraphics[width=3.5in]{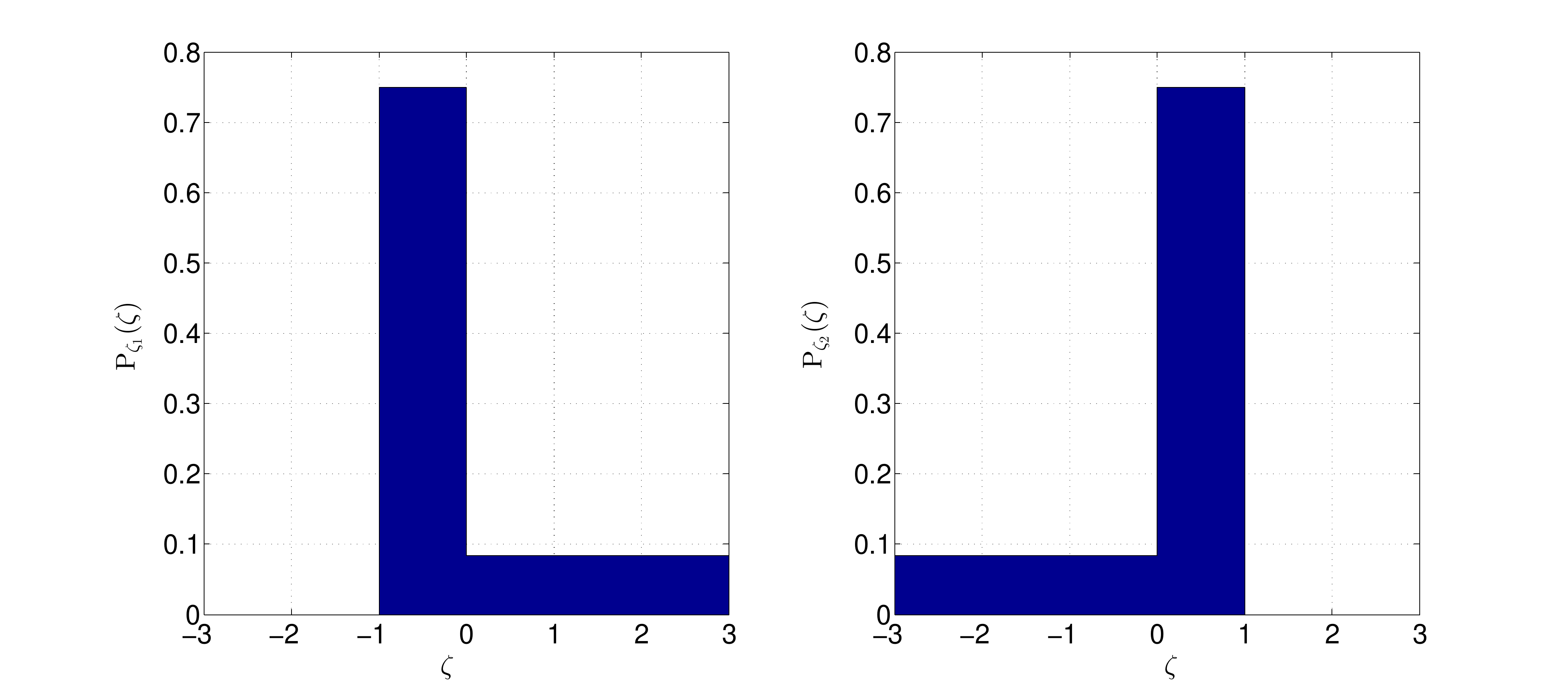}
\caption{True estimate $\mathbb{P}_{\zeta_1}$ and biased estimate $\mathbb{P}_{\zeta_2}$}
\vspace{-.3cm}
\label{fig:numerical}
\end{figure}

We use a simple example to illustrate the impact of an inaccurate mode estimate. We assume random variable $\zeta$ follows distribution $\mathbb{P}_{\zeta_1}$. $\mathbb{P}_{\zeta_2}$ is a biased estimate of $\mathbb{P}_{\zeta_1}$ due to sampling errors. Both distributions are illustrated in Fig.~\ref{fig:numerical}, where each has zero mean and unit variance. However, $\mathbb{P}_{\zeta_1}$ is right-skewed with mode at $-1$ and $\mathbb{P}_{\zeta_2}$ is left-skewed with mode at $1$. Suppose that we try to reformulate $\mathbb{P}_{\zeta}(\zeta\leq z)\geq 90\%$. Based on the given distributions, we find $z\geq 1.8$ from the correct distribution $\mathbb{P}_{\zeta_1}$ and $z\geq 0.925$ from the biased distribution $\mathbb{P}_{\zeta_2}$. In this example, we observe that a misspecified mode estimate could shrink the $90\%$ confidence bound by almost a half and significantly decrease the reliability of the solution to the chance constraint.

\section{Main Results}\label{sec:main_drcc}

\subsection{Assumptions and prior results}
To compute the exact reformulation of distributionally robust chance constraints with various ambiguity sets, we make the following assumptions.
\vspace{.2cm}
\begin{assume}{}\label{asmp: sdp}
For $\mathcal{D}_{\xi}^2$, we assume that
\begin{equation*}
\left(\frac{\alpha+2}{\alpha}\right)(\Sigma-\mu\mu^{\top})\succ\frac{1}{\alpha^2}(\mu-m_{t})(\mu-m_{t})^{\top}.
\end{equation*}
Similarly, for $\mathcal{D}_{\xi}^3$, we assume that, $\forall m\in\Xi,$
\begin{equation*}
\left(\frac{\alpha+2}{\alpha}\right)(\Sigma-\mu\mu^{\top})\succ\frac{1}{\alpha^2}(\mu-m)(\mu-m)^{\top}.
\end{equation*}
\end{assume}
\begin{assume}{}\label{asmp: mode}
For $\mathcal{D}_{\xi}^2$, we assume that $a(x)^{\top}m_{t}\leq b(x)$. Similarly, for $\mathcal{D}_{\xi}^3$, we assume that $a(x)^{\top}m\leq b(x)$, $\forall m\in\Xi$.
\end{assume}
\vspace{.2cm}
Both assumptions are standard in the related literature \cite{hanasusanto2015decision,van2015generalized,van2015distributionally,unimodalBL}. Assumption \ref{asmp: sdp} ensures that the corresponding $\mathcal{D}_{\xi}^i\neq\emptyset$. Assumption \ref{asmp: mode} ensures that the constraint is satisfied at the mode. Furthermore, we assume $\epsilon<0.5$ and $\alpha\geq 1$, since in practice the uncertainties will at least be univariate-unimodal.

Reformulations of \eqref{eq: drcc} under $\mathcal{D}_{\xi}^1$ and $\mathcal{D}_{\xi}^2$ are derived in  previous work.

\begin{thm}{(Theorem 2.2 in \cite{wagner2008stochastic})}\label{thm: drcc_set1}
With $\mathcal{D}_{\xi}^1$, \eqref{eq: drcc} can be exactly reformulated as
\begin{equation}
\sqrt{\left(\frac{1-\epsilon}{\epsilon}\right)a(x)^{\top}(\Sigma-\mu\mu^{\top})a(x)}\leq b(x)-a(x)^{\top}\mu. \label{eq: drcc_set1}
\end{equation}
\end{thm}

\begin{thm}{(Theorem 1 in \cite{unimodalBL})}\label{thm: drcc_set2}
With $\mathcal{D}_{\xi}^2$, \eqref{eq: drcc} can be exactly reformulated as
\begin{align}
&\sqrt{\frac{1-\epsilon-\tau^{-\alpha}}{\epsilon}}\|\Lambda_{t} a(x)\|\leq \tau\left(b(x)-\mu^{\top}a(x)\right)\nonumber\\
&+\left(\tau-\frac{\alpha+1}{\alpha}\right)(\mu-m_{t})^{\top}a(x),\quad\forall\tau\geq\left(\frac{1}{1-\epsilon}\right)^{1/\alpha},\label{eq: drcc_set2}
\end{align}
where $\Lambda_t:=\left(\left(\frac{\alpha+2}{\alpha}\right)(\Sigma-\mu\mu^{\top})-\frac{1}{\alpha^2}(\mu-m_{t})(\mu-m_{t})^{\top}\right)^{1/2}$.
\end{thm}
\vspace{.2cm}
Since parameter $\tau$ has an infinite number of choices, the reformulation in Theorem~\ref{thm: drcc_set2} also involves an infinite number of SOC constraints. Here we obtain a similar result for the generalized ambiguity set $\mathcal{D}^3_{\xi}$.

\subsection{Reformulation for $\mathcal{D}_{\xi}^3$}\label{sec: drcc_set3_1}
We now present the reformation with $\mathcal{D}_{\xi}^3$, which is based on Theorem~\ref{thm: drcc_set2}:
\begin{align}
&\sqrt{\frac{1-\epsilon-\tau^{-\alpha}}{\epsilon}}\|\Lambda a(x)\|\leq \left(\tau-\frac{\alpha+1}{\alpha}\right)(\mu-m)^{\top}a(x)\nonumber\\
&+\tau\left(b(x)-\mu^{\top}a(x)\right),\ \forall\tau\geq\left(\frac{1}{1-\epsilon}\right)^{1/\alpha},\ \forall m\in\Xi,\label{eq: drcc_set3}\\
&\hspace{1.6cm} a(x)^{\top}m\leq b(x),\ \forall m\in\Xi,\label{eq: drcc_asmp2}
\end{align}
where $\Lambda:=\left(\left(\frac{\alpha+2}{\alpha}\right)(\Sigma-\mu\mu^{\top})-\frac{1}{\alpha^2}(\mu-m)(\mu-m)^{\top}\right)^{1/2}$ and \eqref{eq: drcc_asmp2} comes from Assumption~\ref{asmp: mode}.
\vspace{.1cm}

Compared to \eqref{eq: drcc_set2}, \eqref{eq: drcc_set3} is more complicated with two parameters $m$ and $\tau$ each with an infinite number of choices. To solve an optimization problem with \eqref{eq: drcc_set3}, we propose an iterative solving algorithm given in Algorithm~1.

\begin{figure}[h]
 \removelatexerror
  \begin{algorithm}[H]
   \caption{Iterative solving algorithm}
   Initialization: $i=1$, $\tau_0=\bigl(\frac{1}{1-\epsilon}\bigr)^{1/\alpha}$, and $m_0=\{\mbox{any singular point in }\Xi\}$\;
   \vspace{0.1cm}
   Iteration $i$:

   Step 1: Solve the reformulated optimization problem with \eqref{eq: drcc_set3} using $\tau_j$ and $m_j$ for all $j=0,\ldots,i-1$ and obtain optimal solution $x_i^*$. All $\tau_j$ and $m_j$ values are collected from previous iterations;

   Step 2: Find worst case $\tau^*$ and $m^*$ that result in the largest violation of \eqref{eq: drcc_set3} under $x_i^*$: {\bf IF} $m^*$ and $\tau^*$ does not exist, {\bf STOP} and {\bf RETURN} $x_i^*$ as optimal solution; {\bf ELSE} {\bf GOTO} Step 3\;

   Step 3: Set $\tau_i=\tau^*$, $m_i=m^*$, and $i=i+1$\;
  \end{algorithm}
  \label{Fig: algorithm1}
\end{figure}
\noindent Note that the reformulated optimization problem in Step 1 contains only SOC constraints.

\subsection{Step 2 of Algorithm 1}\label{sec: drcc_set3_2}
The challenge is how to efficiently perform Step 2 of Algorithm 1. In the following, we assume $a(x^*_i)\neq0$, otherwise \eqref{eq: drcc_set3} is satisfied with $x^*_i$ regardless of the values of $\tau$ and $m$. Next, we define the following terms
\begin{align*}
&h=a(x^*_i)^{\top}(\mu-m)/\alpha,\  \tilde{c}=b(x^*_i)-\mu^{\top}a(x^*_i),\\
&\tilde{R}=\sqrt{a(x^*_i)^{\top}\left(\frac{\alpha+2}{\alpha}\right)(\Sigma-\mu\mu^{\top})a(x^*_i)}, \\
&g(\tau)=\sqrt{\frac{1-\epsilon-\tau^{-\alpha}}{\epsilon}},\ f(\tau)=-(\alpha\tau-\alpha-1).
\end{align*}
Since $m\in\Xi$, we have $h\in[\underline{h},\overline{h}]$ where
\begin{equation}
\overline{h}=\max_{m\in\Xi}a(x^*_i)^{\top}(\mu-m)/\alpha,\ \underline{h}=\min_{m\in\Xi}a(x^*_i)^{\top}(\mu-m)/\alpha. \label{eq: DRCC_hbar}
\end{equation}
From Assumption~\ref{asmp: sdp}, we have $[\underline{h},\overline{h}]\in(-\tilde{R},\tilde{R})$ and transform \eqref{eq: drcc_set3} into
\begin{align}
&\left[g(\tau)\sqrt{\tilde{R}^2-h^2}+f(\tau)h\right]-\tilde{c}\tau\leq0,\ \forall h \in [\underline{h}, \overline{h}],\ \forall \tau \geq \tau_0.\label{eq: cc-mode}
\end{align}
Since the left side of \eqref{eq: cc-mode} is not jointly convex or concave in $h$ and $\tau$ (see a proof in Appendix~\ref{append1}), we can not find the global maximum value for the left side by simply checking the boundary values or stationary points. Therefore, we propose the following algorithm to efficiently find the global maximum.

We notice that for given a $\tau$ and if $h\in[-\tilde{R},\tilde{R}]$, the maximum value of $g(\tau)\sqrt{\tilde{R}^2-h^2}+f(\tau)h$ equals $\tilde{R}\sqrt{{g(\tau)}^2+{f(\tau)}^2}$ with maximizer $\hat{h}(\tau)=\frac{f(\tau)}{\sqrt{{g(\tau)}^2+{f(\tau)}^2}}\tilde{R}$.
Next, by taking the derivative of $\hat{h}$, we observe that $\hat{h}$ is a strictly decreasing function of $\tau$. Hence, we can compute $\underline{\tau}$ and $\overline{\tau}$ that cause $h$ to reach its boundary values by solving $\hat{h}(\underline{\tau})=\overline{h}\mbox{  and  } \hat{h}(\overline{\tau})=\underline{h}$.
Since $\epsilon< 1/2$ and $\alpha\geq1$, we have $\tau_0<(\alpha+1)/\alpha$ and hence $\hat{h}(\tau_0)=\tilde{R}$. Then we know $[\underline{\tau},\overline{\tau}]>\tau_0$ as $[\underline{h},\overline{h}]< \tilde{R}$.

To efficiently solve these two equalities, we will use a golden section search by first solving for $\overline{\tau}$ on $[\tau_0,\infty]$ and then for $\underline{\tau}$ on $[\tau_0,\overline{\tau}]$. To efficiently apply a golden section search on $\overline{\tau}$, we need to find a finite upper bound instead of $\infty$. The following lemma describes the selection of the finite upper bound $\tau_1$ and the best region to conduct the golden section search.
\vspace{0.2cm}
\begin{lemma}{}\label{lem: DRCC1}
If $\underline{h}\geq0$, $\tau_1=\frac{\alpha+1}{\alpha}$. The golden section search of $\overline{\tau}$ can be conducted on $[\tau_0,\frac{\alpha+1}{\alpha}]$. If $\underline{h}<0$, $\tau_1=-\left(\underline{h}\sqrt{\frac{1-\epsilon}{\epsilon(\tilde{R}^2-\underline{h}^2)}}-(\alpha+1)\right)/\alpha$. The search can be conducted on $[\frac{\alpha+1}{\alpha},\tau_1]$. The proof is given in Appendix~\ref{append2}.
\end{lemma}
\vspace{0.2cm}
Furthermore, from Assumption \ref{asmp: mode}, we have $\tilde{c}\geq -\alpha\underline{h}\geq-\alpha\overline{h}$.

Based on the threshold values $\overline{\tau}$ and $\underline{\tau}$, we divide our discussion into three cases.
\vspace{0.2cm}

\noindent\textbf{Case 1}: If $\tau\in[\tau_0,\underline{\tau}]$, $h^*=\overline{h}$. Then from \eqref{eq: cc-mode}, we find
$$g(\tau)\sqrt{\tilde{R}^2-\overline{h}^2}+f(\tau)\overline{h}-\tilde{c}\tau\leq0.$$
Then, we transform the above constraint into the following equivalent form:
\begin{equation}
F_1(\tau)=C_1g(\tau)-(\tilde{c}+\alpha\overline{h})\tau+(\alpha+1)\overline{h}\leq 0,\label{cc-mode-3}
\end{equation}
where $C_1=\sqrt{\tilde{R}^2-\overline{h}^2}\geq0$. The left side of \eqref{cc-mode-3} is concave on $\tau$. Define the derivative of the left side as $F_1^{\prime}(\tau)=C_1g^{\prime}(\tau)-(\tilde{c}+\alpha\overline{h}).$
We observe that $F_1^{\prime}(\tau_0)>0$ as $g^{\prime}(\tau_0)\to\infty$ and
\begin{enumerate}
  \item if $F_1^{\prime}(\underline{\tau})\leq 0$, $\tau^*$ is the unique solution of $F_1^{\prime}(\tau)=0$ within the domain $[\tau_0,\underline{\tau}]$;
  \item else if $F_1^{\prime}(\underline{\tau})> 0$, $\tau^*=\underline{\tau}$.
\end{enumerate}
\vspace{0.2cm}

\noindent\textbf{Case 2}: If $\tau\in[\underline{\tau},\overline{\tau}]$, $h^*=\hat{h}(\tau)$. Then from \eqref{eq: cc-mode}, we find
$$\tilde{R}\sqrt{{g(\tau)}^2+{f(\tau)}^2}-\tilde{c}\tau\leq0.$$
The above problem is a one-dimensional problem on $\tau$. We transform it into the following form:
\begin{equation}
F_2(\tau)=\tilde{R}^2\left({g(\tau)}^2+{f(\tau)}^2\right)-\tilde{c}^2\tau^2\leq 0.\label{cc-mode-4}
\end{equation}
We observe that $F_2(\tau)$ is differentiable on $[\underline{\tau},\overline{\tau}]$. Then, we know that the extreme value of $F_2(\tau)$ happens at the critical points (boundary points $\underline{\tau}$, $\overline{\tau}$ or $\tau_i$ such that that $F_2^{\prime}(\tau_i)=0$). In the following numerical analysis, we present efficient ways to find $\tau^*$ which maximize the left side of \eqref{cc-mode-4}.

The first and second derivative of the left side of \eqref{cc-mode-4} are
\begin{align*}
F_2^{\prime}(\tau)&=\tilde{R}^2\left(\frac{\alpha}{\epsilon}\tau^{-\alpha-1}+2\alpha(\alpha\tau-\alpha-1)\right)-2\tilde{c}^2\tau\\
&=\frac{\alpha\tilde{R}^2}{\epsilon}\tau^{-\alpha-1}+(2\alpha^2\tilde{R}^2-2\tilde{c^2})\tau-2\tilde{R}^2\alpha(\alpha+1),\\
F_2^{\prime\prime}(\tau)&=-\frac{\alpha\tilde{R}^2(\alpha+1)}{\epsilon}\tau^{-\alpha-2}+(2\alpha^2\tilde{R}^2-2\tilde{c^2}).
\end{align*}
Given this, there are two conditions.

\textit{Condition 1:} If $2\alpha^2\tilde{R}^2-2\tilde{c^2}\leq 0$, $F_2^{\prime}(\tau)$ is monotonically decreasing on $\tau$ and $F_2(\tau)$ is concave on $\tau$. Then,
\begin{enumerate}
\item if $F_2^{\prime}(\underline{\tau})\leq 0$, $\tau^*=\underline{\tau}$;
\item  else if $F_2^{\prime}(\underline{\tau})>0$ and $F_2^{\prime}(\overline{\tau})\leq 0$, $\tau^*$ is the unique solution of $F_2^{\prime}(\tau)=0$ within the domain $[\underline{\tau},\overline{\tau}]$.
\item else if  $F_2^{\prime}(\overline{\tau})>0$, $\tau^*=\overline{\tau}$.
\end{enumerate}

\textit{Condition 2:} If $2\alpha^2\tilde{R}^2-2\tilde{c^2}> 0$, $F_2^{\prime\prime}(\tau)$ is monotonically increasing on $\tau$ and $F_2^{\prime}(\tau)$ is convex on $\tau$. Then,
\begin{enumerate}
\item  if $F_2^{\prime\prime}(\overline{\tau})\leq 0$, $F_2^{\prime}(\tau)$ is decreasing within the domain. To find $\tau^*$, we follow the same discussions as in Condition~1;
\item  else if $F_2^{\prime\prime}(\overline{\tau})>0$ and $F_2^{\prime\prime}(\underline{\tau})\leq 0$, $F_2^{\prime}(\tau)$ is first decreasing and then increasing. Define $F_s=F_2^{\prime}(\tau_s)$ where $F_2^{\prime\prime}(\tau_s)=0$ within the domain $[\underline{\tau},\overline{\tau}]$, $F_l=F_2^{\prime}(\underline{\tau})$, and $F_u=F_2^{\prime}(\overline{\tau})$. Then,

\begin{enumerate}
  \item If $0\leq F_s$, $\tau^*=\overline{\tau}$.
  \item If $F_s\leq 0 \leq F_l \leq F_u$ or $F_s\leq 0 \leq F_u \leq F_l$, $\tau*=\overline{\tau}$ or the unique solution of $F_2^{\prime}(\tau)=0$ within the domain $[\underline{\tau},\tau_s]$ that maximizes $F_2(\tau)$.
  \item If $F_s\leq F_l\leq 0\leq F_u$, $\tau^*=\underline{\tau}$ or $\overline{\tau}$ that maximizes $F_2(\tau)$.
  \item If $F_s\leq F_u \leq 0 \leq F_l$, $\tau^*$ equals the unique solution of $F^{\prime}(\tau)=0$ within the domain $[\underline{\tau},\tau_s]$.
  \item  If $F_s\leq F_l\leq F_u\leq 0$ or $F_s\leq F_u\leq F_l\leq 0$, $\tau^*=\underline{\tau}$.
\end{enumerate}

\item else if $F_2^{\prime\prime}(\underline{\tau})> 0$, $F_2(\tau)$ is convex on $\tau$. $\tau^*=\underline{\tau}$ or $\overline{\tau}$ that maximizes $F_2(\tau)$.
\end{enumerate}
\vspace{0.2cm}

\noindent\textbf{Case 3}: If $\tau\in[\overline{\tau},\infty]$, $h^*=\underline{h}$. Then from \eqref{eq: cc-mode}, we find
\begin{equation}
g(\tau)\sqrt{\tilde{R}^2-\underline{h}^2}+f(\tau)\underline{h}-\tilde{c}\tau\leq0,\nonumber
\end{equation}
which we transform into the following equivalent form
\begin{equation}
F_3(\tau)=C_3g(\tau)-(\tilde{c}+\alpha\underline{h})\tau+(\alpha+1)\underline{h}\leq 0,\label{cc-mode-5}
\end{equation}
where $C_3=\sqrt{\tilde{R}^2-\underline{h}^2}$. Define the derivative of the left hand side of \eqref{cc-mode-5} as $F_3^{\prime}(\tau)=C_3g^{\prime}(\tau)-(\tilde{c}+\alpha\underline{h}).$
Then $F_3(\tau)$ is concave on $\tau$ and as $\tau\to\infty$, $F_3^{\prime}(\tau)\leq0$. Then,
\begin{enumerate}
\item if $\tilde{c}+\alpha\underline{h}=0$, $F_3(\tau)$ is an increasing function and $\tau^*=\infty$;
\item else if $\tilde{c}+\alpha\underline{h}>0$, as $\tau\to\infty$, $F_3^{\prime}(\tau)<0$. Based on the concavity of $F_3(\tau)$, we find
\begin{enumerate}
  \item if $F_3^{\prime}(\overline{\tau})\leq 0$, $\tau^*=\overline{\tau}$;
  \item else if $F_3^{\prime}(\overline{\tau})> 0$, $\tau^*$ equals the unique solution of $F_3^{\prime}(\tau)=0$ within the domain $[\overline{\tau},\infty]$.
\end{enumerate}
\end{enumerate}

To efficiently apply the golden section search, we determine an effective finite upper bound instead of $\infty$. Let the effective upper bound be $\tau_2$, we have
$$F_3^{\prime}(\tau_2)=C_3g^{\prime}(\tau_2)-(\tilde{c}+\alpha\underline{h})\leq0.$$
\vspace{0.2cm}
\begin{lemma}{}\label{lem: DRCC2}
A feasible selection of $\tau_2$ is
$$\tau_2=\left[\frac{-1+\sqrt{1+4(1-\epsilon)C_2}}{2C_2}\right]^{-\frac{1}{\alpha}},$$
where $C_2=\frac{\alpha^2C_3^2}{4\epsilon (\tilde{c}+\alpha\underline{h})^2}$. The proof is given in Appendix~\ref{append3}.
\end{lemma}
\vspace{0.2cm}
Then, instead of a search on $[\overline{\tau},\infty]$, we only need to search on $[\overline{\tau},\tau_2]$.

Combining all three cases, we can find the overall worst case $\tau^*$ and $h^*$ given $x^*_i$. If \eqref{eq: cc-mode} is satisfied with these parameters, then there is no violated constraint in Step 2 of Algorithm 1. If \eqref{eq: cc-mode} is not satisfied, we need to use the worst case $\tau^*$ and $m^*$ in Step 3 and the iteration continues. Depending on how we define $\Xi$, $m^*$ are different functions of $h^*$.

\subsection{Candidates of $\Xi$}\label{sec: Xi}
In this section, we demonstrate how the selection of $\Xi$ affects the determination of $\underline{h}$, $\overline{h}$, and $m^*$. Specifically, we give two examples of $\Xi$ and show how to exactly reformulate \eqref{eq: drcc_asmp2} (i.e., Assumption~\ref{asmp: mode}) and how to calculate $\underline{h}$ and $\overline{h}$, given $x^*_i$. Furthermore, we show how to find the worst case $m^*$ from $h^*$.
\vspace{0.2cm}

\noindent\textbf{Rectangular Support}: We assume that $\Xi=[\underline{k},\overline{k}]$ and hence we can reformulate \eqref{eq: drcc_asmp2} as
\begin{equation}
a(x)^{\top}\left(\frac{\underline{k}+\overline{k}}{2}\right)+\left|a(x)\right|^{\top}\left(\frac{\overline{k}-\underline{k}}{2}\right)\leq b(x).\label{eq: DRCC_rec_asmp2}
\end{equation}
Furthermore, given $x^*_i$, we have the following relationships due to \eqref{eq: DRCC_hbar}.
\begin{align}
&\underline{h}=\left[a(x^*_i)^{\top}\left(\mu-\frac{\underline{k}+\overline{k}}{2}\right)-|a(x^*_i)|^{\top}\left(\frac{\overline{k}-\underline{k}}{2}\right)\right]/\alpha,\label{eq: DRCC_rec1}\\
&\overline{h}=\left[a(x^*_i)^{\top}\left(\mu-\frac{\underline{k}+\overline{k}}{2}\right)+|a(x^*_i)|^{\top}\left(\frac{\overline{k}-\underline{k}}{2}\right)\right]/\alpha. \label{eq: DRCC_rec2}
\end{align}
Based on \eqref{eq: DRCC_rec1} and \eqref{eq: DRCC_rec2}, if we have the worst case $h^*$, we find the worst case $m^*$ by solving \eqref{eq: DRCC_rec3} for $\lambda_r$ and substituting in \eqref{eq: DRCC_rec4}:
\begin{align}
&h^*=\left[a(x^*_i)^{\top}\left(\mu-\frac{\underline{k}+\overline{k}}{2}\right)+\lambda_r|a(x^*_i)|^{\top}\left(\frac{\overline{k}-\underline{k}}{2}\right)\right]/\alpha,\label{eq: DRCC_rec3}\\
&m^*=\left(\frac{\underline{k}+\overline{k}}{2}\right)-\lambda_r\textbf{sign}\left(a(x^*_i)\right)\left(\frac{\overline{k}-\underline{k}}{2}\right),\label{eq: DRCC_rec4}
\end{align}
where $\textbf{sign}\left(a(x)\right)$ returns a diagonal matrix whose diagonal elements equal the sign of each elements in $a(x)$.
\vspace{0.2cm}

\noindent\textbf{Ellipsoidal Support}: We assume that $\Xi = \{m: m = m_c + P^{1/2}u, ||u||_2 \leq 1\}$, where $P\succ0$. Then we can reformulate \eqref{eq: drcc_asmp2} as
\begin{equation}
a(x)^{\top}m_c+\left\|P^{1/2}a(x)\right\|_2\leq b(x).\label{eq: DRCC_ellip_asmp2}
\end{equation}
Furthermore, due to \eqref{eq: DRCC_hbar}, we have the following relationships:
\begin{align}
&\underline{h}=\left[a(x^*_i)^{\top}(\mu-m_c)-\left\|P^{1/2}a(x^*_i)\right\|_2\right]/\alpha,\label{eq: DRCC_ellip1}\\
&\overline{h}=\left[a(x^*_i)^{\top}(\mu-m_c)+\left\|P^{1/2}a(x^*_i)\right\|_2\right]/\alpha. \label{eq: DRCC_ellip2}
\end{align}
Next, if we have the worst case $h^*$, we find the worst case $m^*$ directly by solving \eqref{eq: DRCC_ellip3} for $\lambda_e$ and substituting in \eqref{eq: DRCC_ellip4}:
\begin{align}
&h^*=\left[a(x^*_i)^{\top}(\mu-m_c)-\lambda_e\left\|P^{1/2}a(x^*_i)\right\|_2\right]/\alpha,\label{eq: DRCC_ellip3}\\
&m^*=m_c+\lambda_e\frac{Pa(x^*_i)}{\left\|P^{1/2}a(x^*_i)\right\|_2}. \label{eq: DRCC_ellip4}
\end{align}

\section{Case Study}\label{sec:casestudy}
\subsection{Simulation Setup}
We consider the DC OPF problem from \cite{unimodalBL}. We assume that the system has two wind power plants with wind forecast error $\tilde{w}=[\tilde{w}_1,\tilde{w}_2]^{\top}$. With $N_G$ generators and $N_B$ buses, the design variables are generation $P_G \in \mathbb R^{N_G}$, up and down reserve capacities $R_G^{up} \in \mathbb R^{N_G}, R_G^{dn} \in \mathbb R^{N_G}$, and a distribution vector $d_G \in \mathbb R^{N_G}$, which determines the real-time reserve provision from each generator used to balance the wind forecast error. The full problem formulation is as follows.
\begin{subequations}
\begin{align}
\min\ &P_G^T[C_1]P_G+C_2^TP_G+C_R^T(R_G^{up}+R_G^{dn})\label{eq:obj}\\
\mathrm{s.t.}
&-P_l\leq AP_\text{inj}\leq P_l\label{cons:line}\\
&R_G=-d_G(\tilde{w}_1+\tilde{w}_2)\label{cons:resv}\\
&P_{\text{inj}}=C_G(P_G+R_G)+C_W(P_W^f+\tilde{w})-C_LP_L\label{cons:inj}\\
&\underline{P}_G\leq P_G+R_G\leq \overline{P}_G\label{cons:gen}\\
&-R_G^{dn}\leq R_G \leq R_G^{up}\label{cons:resvcap}\\
&\mathbf{1}_{1 \times N_G} d_G=1\label{cons:dist}\\
&\mathbf{1}_{1 \times N_B}(C_GP_G+C_WP_W^f-C_LP_L)=0\label{cons:bal}\\
&P_{G}\geq \mathbf{0}_{N_G \times 1},\ d_G\geq \mathbf{0}_{ N_G \times 1} \label{cons:nneg1}\\
& R_G^{up}\geq \mathbf{0}_{N_G \times  1}, \ R_G^{dn}\geq \mathbf{0}_{N_G \times 1} \label{cons:nneg2}
\end{align}
\end{subequations}
where $[C_1] \in \mathbb R^{N_G \times N_G}$, $C_2 \in \mathbb R^{N_G}$, and $C_R \in \mathbb R^{N_G}$ are cost parameters. Constraint \eqref{cons:line} bounds the power flow, which is calculated from the power injections $P_\text{inj}$ defined in \eqref{cons:inj}  and the parameter matrix $A$, by the line limits $P_l$. Constraint \eqref{cons:resv} computes the real-time reserve usage $R_G$ for each generator. In \eqref{cons:inj} $P_W^f$ is the wind forecast, $P_L$ is the load, and $C_G$, $C_W$, and $C_L$ are  matrices that map generators, wind power plants, and loads to buses; \eqref{cons:gen} restricts generation  to within its limits $[\underline{P}_G, \overline{P}_G]$; \eqref{cons:resvcap} restricts $R_G$ by the reserve capacity; \eqref{cons:dist}, \eqref{cons:bal} enforce power balance with and without wind forecast error; and \eqref{cons:nneg1}, \eqref{cons:nneg2} ensure all decision variables are non-negative.

We test our approach on a modified IEEE 30-bus system with network and cost parameters from \cite{matpower}. We set $C_R=10C_2$. We add the wind power plants to buses 22 and 5 and set $P_W^f= [66.8, 68.1]$ MW. We use the same wind power forecast uncertainty data ($10000$ scenarios) as in Section~\ref{sec: winderror}. We congest the system by increasing each load by $50\%$ and reducing the limit of the line connecting buses 1 and 2 to 30 MW. All optimization problems are solved using CVX with the Mosek solver \cite{cvx1,cvx2}.

To construct the ambiguity sets, unlike in Section~\ref{sec: winderror}, the outliers are used when estimating the statistical parameters (first moment $\mu$, second moment $\Sigma$, and the set of the mode $\Xi$) and evaluating the reliability of the solution. We set $\epsilon=5\%$, $\alpha=1$, and assume $\Xi$ is a rectangular set.

\subsection{Additional Ambiguity Sets}
We benchmark against two additional ambiguity sets from related work.\\
 {\bf Ambiguity set 4:} (moment and unimodality with fixed mode at the mean \cite{stellato2014data})
\begin{align}
&\mathcal{D}_{\xi}^4 := \bigl\{\mathbb{P}_{\xi} \in \mathcal{P}^n_{\alpha} \cap \mathcal{D}_{\xi}^1:\ \mathcal{M}(\xi)=\mu\bigr\}. \label{Dset4}
\end{align}
{\bf Ambiguity set 5:} (moment and unimodality with $\alpha=1$ and arbitrary mode \cite{roald1})
\begin{align}
&\mathcal{D}_{\xi}^5 := \left\{\mathbb{P}_{\xi} \in \mathcal{P}^n_{1} \cap \mathcal{D}_{\xi}^1 \right\}. \label{Dset5}
\end{align}
Set $\mathcal{D}_{\xi}^4$ is a special case of $\mathcal{D}_{\xi}^2$ with the mode at the mean, while $\mathcal{D}_{\xi}^5$ is a special case of $\mathcal{D}_{\xi}^3$ with $\alpha=1$ and $\Xi:=\{\mbox{all possible values of }\mathcal{M}(\xi)\}$ that is an ellipsoidal set based on $\mu$ and $\Sigma$ as shown in Assumption~\ref{asmp: sdp}. In other words, our $\mathcal{D}_{\xi}^3$ is more general than $\mathcal{D}_{\xi}^{2}$, $\mathcal{D}_{\xi}^{4}$, and $\mathcal{D}_{\xi}^{5}$. The reformulations of  $\mathcal{D}_{\xi}^4$ and $\mathcal{D}_{\xi}^5$ are simpler than $\mathcal{D}_{\xi}^3$ with a single SOC constraint
\begin{equation}
\mathcal{K}\sqrt{a(x)^{\top}(\Sigma-\mu\mu^{\top})a(x)}\leq b(x)-a(x)^{\top}\mu, \label{eq: drcc_set45}
\end{equation}
where $\mathcal{K}$ can be found in \cite{stellato2014data} for $\mathcal{D}_{\xi}^4$ and in \cite{roald1} for $\mathcal{D}_{\xi}^5$.

\subsection{Simulation Results}
\subsubsection{Estimation of $\Xi$}\label{sec: support}
\begin{figure}
\centering
\vspace{.3cm}
\includegraphics[width=3.5in]{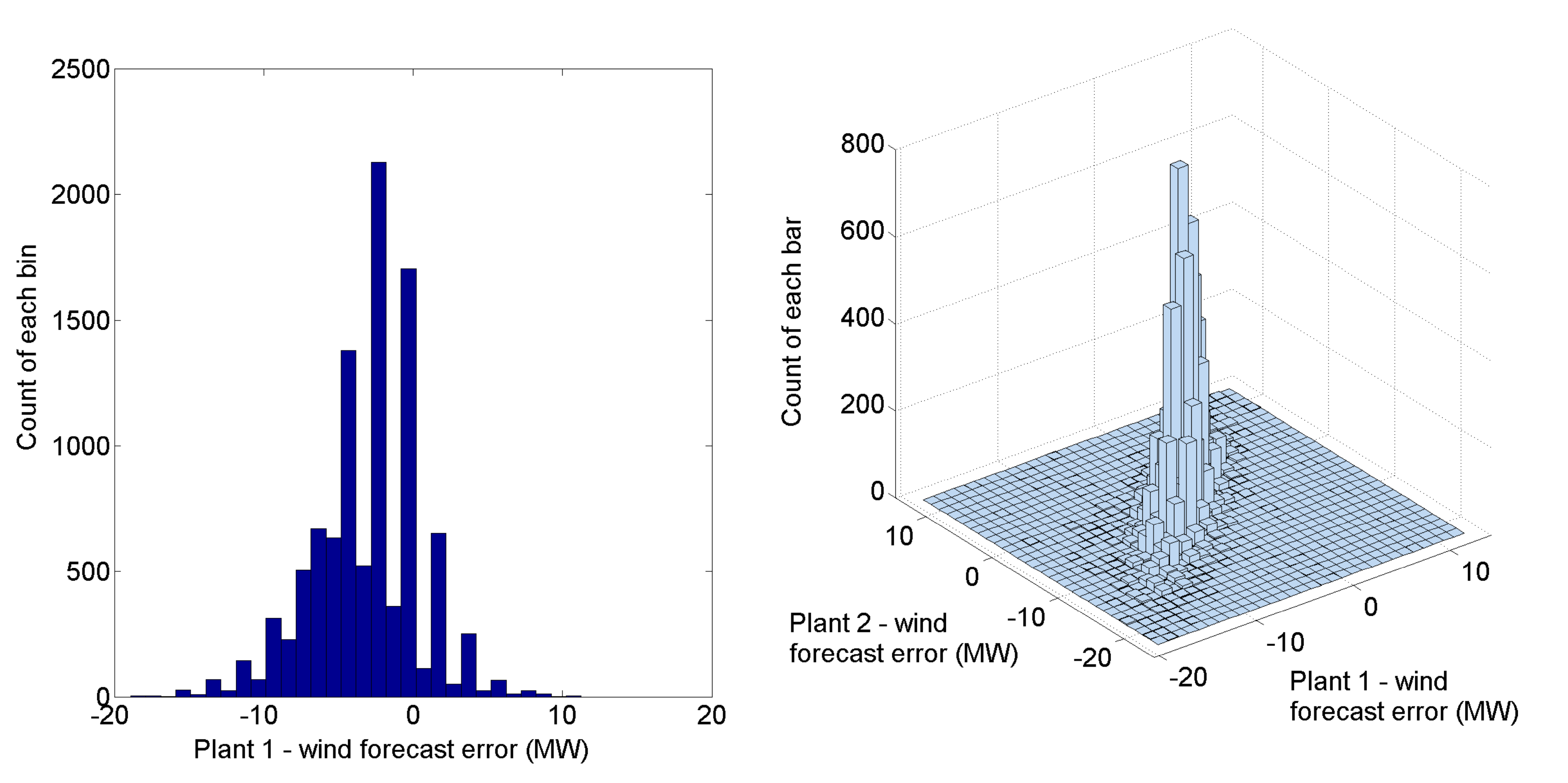}
\caption{Histogram of univariate and bivariate wind forecast errors ($30$ bins).}
\vspace{-.3cm}
\label{fig:windunimodal2}
\end{figure}

We next analyze how  the data size of each sample $N_{data}$ and the number of bins within the histogram $N_{bin}$ affect the estimate of the mode support.  Figure~\ref{fig:windunimodal2} shows that if we change $N_{bin}$ from 15 to 30 the histograms no longer show a unimodal distribution, as compared to Fig.~\ref{fig:windunimodal}. The problem is exacerbated as $N_{bin}$ grows.

\begin{figure}
\centering
\vspace{.3cm}
\includegraphics[width=3.5in]{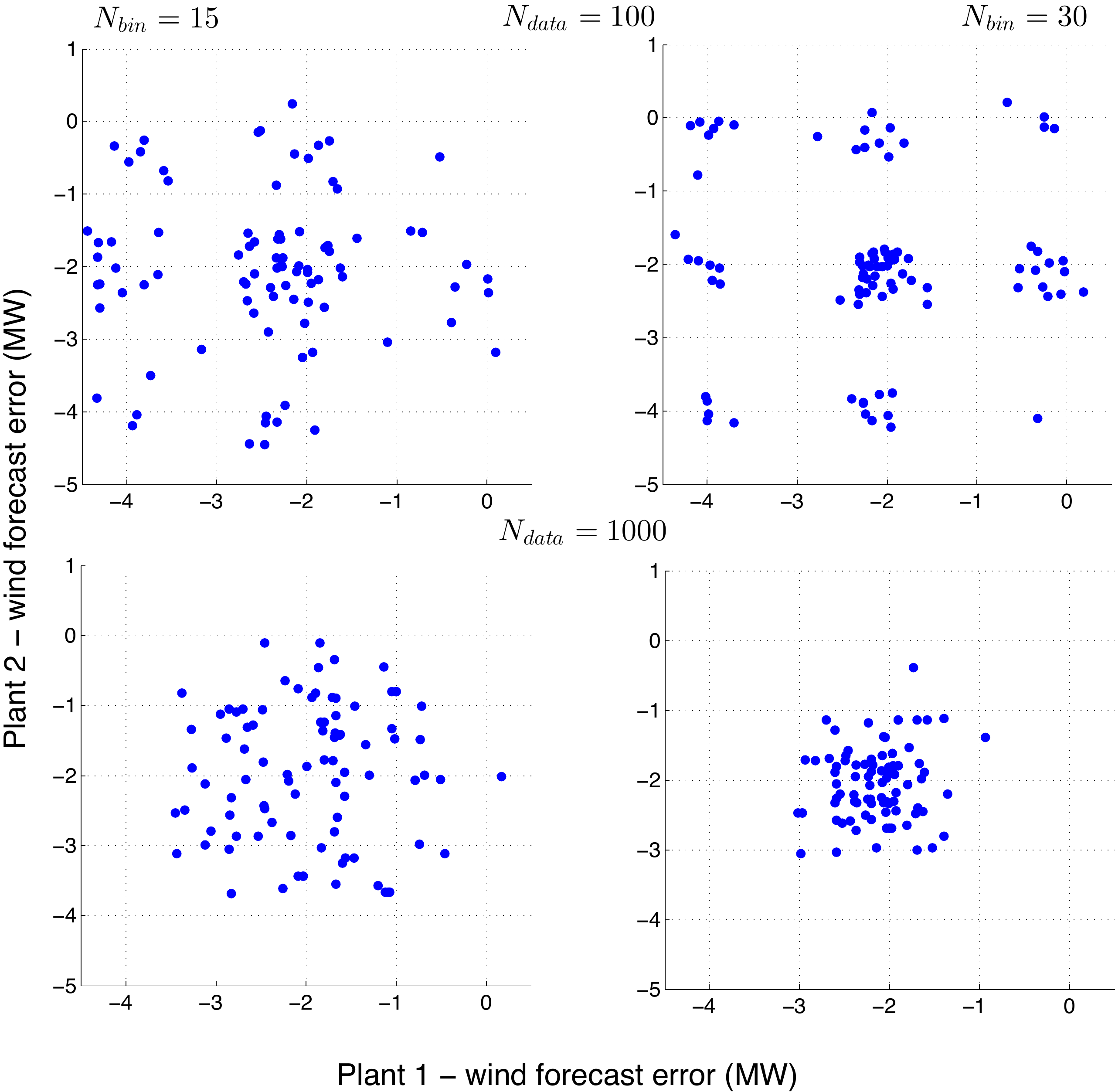}
\caption{Mode values from samples with different $N_{data}$ and $N_{bin}$. Data is sampled from the full data pool.}
\vspace{-.3cm}
\label{fig:windscatter2}
\end{figure}

We next explore the impact of the size of the data pool. We first use the entire data pool to select 100 samples with different data sizes $N_{data}$ ($100$ and $1000$) and number of bins $N_{bin}$ ($15$ and $30$) and show scatter plots of the mode values in Fig.~\ref{fig:windscatter2}. As $N_{data}$ gets larger, the mode values are more condensed and hence more accurate. When $N_{bin}=30$ and $N_{data}=100$ mode values appear in several disjoint regions, but this disjointness is mitigated as $N_{data}$ increases to $1000$. Based on the scatter plots, we determined the parameters $\underline{k},\overline{k}$ of the four rectangular sets $\Xi$ used in $\mathcal{D}_{\xi}^{3}$. The results are given in Table~\ref{tab: recparam}.

\begin{table}
\centering
\caption{Full pool: $\underline{k}$ and $\overline{k}$ (MW) of Four Rectangular Sets $\Xi$. }
\label{tab: recparam}
\begin{tabular}{cc|cc|cc}
\hline
                    &  & \multicolumn{2}{c|}{$N_{data}=100$}  & \multicolumn{2}{c}{$N_{data}=1000$} \\
$N_{bin} $         &             & $\underline{k}$         & $\overline{k}$                & $\underline{k}$          & $\overline{k}$         \\ \hline\hline
\multirow{2}{*}{15} & Plant 1     &-4.44            &0.10                       &-3.45             &   0.17         \\
                    & Plant 2       & -4.45           & 0.24                    & -3.69            &   -0.11         \\ \hline
\multirow{2}{*}{30} & Plant 1      &   -4.36         &  0.19                  &   -3.02          &      -0.93      \\
                    & Plant 2      &   -4.22         &  0.22                    &   -3.06          &  -0.39          \\ \hline
\end{tabular}

\vspace{.5cm}
\caption{Partial pool: $\underline{k}$ and $\overline{k}$ (MW) of Four Rectangular Sets $\Xi$}
\begin{tabular}{cc|cc|cc}
\hline
                    &  & \multicolumn{2}{c|}{$N_{data}=50$}  & \multicolumn{2}{c}{$N_{data}=200$} \\
$N_{bin} $         &             & $\underline{k}$         & $\overline{k}$                & $\underline{k}$          & $\overline{k}$         \\ \hline\hline
\multirow{2}{*}{10} & Plant 1     &-4.77            &0.58                       &-3.52             &   0.09         \\
                    & Plant 2      & -5.05           & 0.44                    & -4.43            &   0.06         \\ \hline
\multirow{2}{*}{20} & Plant 1       &   -5.82         &  0.06                  &   -4.36          &      -0.09      \\
                    & Plant 2      &   -5.76         &  0.04                    &   -3.86          &  0.19          \\ \hline
\end{tabular}
\label{tab: recparam2}
\end{table}

\begin{figure}
\centering
\vspace{.3cm}
\includegraphics[width=3.4in]{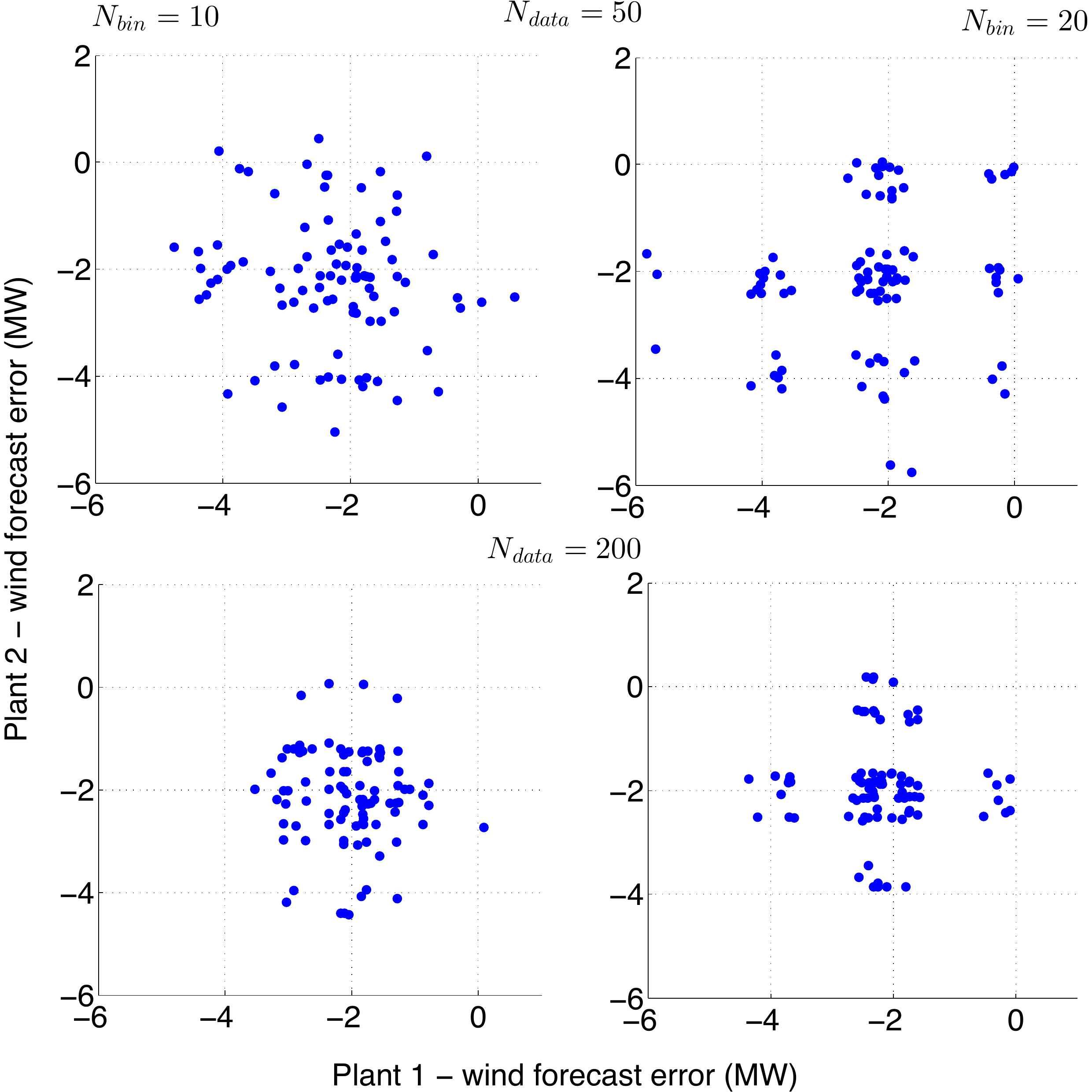}
\caption{Mode values from samples with different $N_{data}$ and $N_{bin}$. Data is sampled from the partial data pool.}
\vspace{-.3cm}
\label{fig:windscatter3}
\end{figure}

We repeated the analysis using only a partial data pool, specifically, we randomly selected 1000 data from the full pool to comprise the partial pool. We also use different choices of $N_{data}$ and $N_{bin}$. The scatter plots are shown in Fig.~\ref{fig:windscatter3} and parameter values for $\Xi$ are given in Table~\ref{tab: recparam2}.

\subsubsection{Objective Costs}\label{sec: cost}
We next analyze the objective costs and the optimal reserve capacities using different ambiguity sets. The results are summarized in Table~\ref{tab:cost}. In all case studies, since we focus on mode misspecification not moment misspecification, moments are calculated using the full or partial data pool and all ambiguity sets use the same moments.

\begin{table*}
\centering
\caption{Objective costs and reserve capacities}
\label{tab:cost}
\begin{tabular}{l|c|cccccc|cccc|c|c}
\hline
       \multirow{2}{*}{Full pool}           & \multirow{2}{*}{$\mathcal{D}_{\xi}^{1}$} & \multicolumn{6}{c|}{$\mathcal{D}_{\xi}^{2}$} & \multicolumn{4}{c|}{$\mathcal{D}_{\xi}^{3}$}&  \multirow{2}{*}{$\mathcal{D}_{\xi}^{4}$} & \multirow{2}{*}{$\mathcal{D}_{\xi}^{5}$} \\
       &&  M$1$ & M$2$ & M$3$ & M$4$& M$5$& M$6$\ & $\Xi_1$ & $\Xi_2$ & $\Xi_3$ & $\Xi_4$ & &  \\ \hline\hline
Total Cost        & 26160 & 19440 & 19546  & 18993 & 19547 & 19526 & 19542 & 19949 &  19818  & 19982 & 19896 &19818  & 19982 \\
Generation Cost    & 13032 & 11515  &  11504 & 11506 & 11481  & 11491 & 11506 & 11522 &  11522  & 11522 &11522 &  11514& 11522 \\
Reserve Cost      & 13129 & 7925  &  8042   & 7488  & 8065  & 8035  & 8036 & 8427  &  8296   & 8460  & 8373 &  8304& 8460 \\
Up Reserve  (MW)   & 38.8  & 26.8   &  26.2   & 26.3  & 25.1  & 26.1  & 26.1 & 27.1  &  26.8   & 27.1  & 27.0&  26.7& 27.1 \\
Down Reserve (MW)  & 26.9  & 12.9   &  14.0   & 11.1  & 15.2   & 14.1   & 14.1 & 15.1  &  14.7   & 15.2  & 14.9 &  14.8& 15.2 \\ \hline
\end{tabular}

\vspace{0.3cm}

\begin{tabular}{l|c|cccccc|cccc|c|c}
\hline
       \multirow{2}{*}{Partial pool}           & \multirow{2}{*}{$\mathcal{D}_{\xi}^{1}$} & \multicolumn{6}{c|}{$\mathcal{D}_{\xi}^{2}$} & \multicolumn{4}{c|}{$\mathcal{D}_{\xi}^{3}$}&  \multirow{2}{*}{$\mathcal{D}_{\xi}^{4}$} & \multirow{2}{*}{$\mathcal{D}_{\xi}^{5}$} \\
       &&  M$1$ & M$2$ & M$3$ & M$4$& M$5$& M$6$\ & $\Xi_5$ & $\Xi_6$ & $\Xi_7$ & $\Xi_8$ & &  \\ \hline\hline
Total Cost        & 21845 & 16759   & 16740  & 15250 & 15915 & 16735 & 16718 & 17883 & 17909  & 17876 & 17909 & 16791    & 17949\\
Generation Cost    & 11713 & 11376   & 11376  & 11324 &11295  & 11336 & 11367 & 11420 &  11420 & 11420 &11420 &  11371 & 11420 \\
Reserve Cost      & 10132 & 5383    &  5364  & 3926  & 4620 & 5399 & 5351 & 6463 &  6489  & 6456  & 6489 &  5420  & 6529\\
Up Reserve (MW)    & 31.1  & 20.1  &  20.1   & 17.5  & 16.1 & 19.4 & 19.1 & 22.0 &  22.0  & 21.9     & 22.0 &  19.2     & 22.1 \\
Down Reserve (MW)  & 19.6 & 6.8     &  6.7    & 2.1   & 7.0 & 7.6 & 7.7 & 10.3   &  10.4  & 10.4    & 10.4 &  7.9       & 10.6 \\ \hline
\end{tabular}
\end{table*}

For ambiguity set $\mathcal{D}_{\xi}^{2}$, we perform tests with the following six fixed mode estimates.
\begin{itemize}
\item M1: mode determined using the full (partial) data pool with histogram of $15$ ($10$) bins. This case demonstrates the performance of $\mathcal{D}_{\xi}^{2}$ with an accurate mode estimate.
\item M2: mode determined using the full (partial) data pool with histogram of $30$ ($20$) bins. This case shows how $N_{bin}$ affects the result.
\item M3-6: combinations of the largest $\overline{k}$ and the smallest $\underline{k}$ of both plants from Table~\ref{tab: recparam} (full pool) and  Table~\ref{tab: recparam2} (partial pool). These cases demonstrate the affect of outlying data samples.
\end{itemize}
\vspace{0.2cm}

For ambiguity set $\mathcal{D}_{\xi}^3$, we perform tests with different $\Xi$, specifically, $\Xi_1:100\times 15, \Xi_2: 1000\times 15, \Xi_3:100\times 30, \Xi_4:1000\times 30, \Xi_5: 50\times 10, \Xi_6: 200\times 10,  \Xi_7: 50\times 20,$ and $\Xi_8: 200\times 20$, where the first number refers to $N_{data}$ and the second number refers to $N_{bin}$. In each case we use the parameters $\underline{k}, \overline{k}$ from Tables~\ref{tab: recparam} and \ref{tab: recparam2}.

As shown in Table~\ref{tab:cost}, $\mathcal{D}_{\xi}^{1}$ has the highest objective cost since it does not include the assumption of unimodality. The cost of $\mathcal{D}_{\xi}^{2}$ varies with the mode estimate. We observe opposite variations on the total up and down reserve capacities since different mode estimates lead to different estimates of the skewness of the uncertainty distribution. Comparing M1 and M2 to M3-6 we see that inaccurate estimation of the mode could lead to either higher or lower costs. Furthermore, results for M1 and M2 are significantly different demonstrating the effect of different choices of $N_{bin}$.

The costs of $\mathcal{D}_{\xi}^{3}$ are higher than those of $\mathcal{D}_{\xi}^{2}$ since the solution is designed to cope with mode misspecification. The costs do not vary significantly as a function of $N_{bin}$ and $N_{data}$. For a given $N_{bin}$, as $N_{data}$ increases, the costs decrease since the mode estimates are more closely clustered.

The cost of  $\mathcal{D}_{\xi}^{4}$ is higher than the costs of $\mathcal{D}_{\xi}^{2}$ with M1, demonstrating the benefit in allowing the mode to be different than the mean. The cost of $\mathcal{D}_{\xi}^{5}$ is close to that of $\mathcal{D}_{\xi}^{3}$ with $\Xi_3$ since the mode estimates are widely distributed in this case; however, the cost of all other $\mathcal{D}_{\xi}^{3}$ is below that of $\mathcal{D}_{\xi}^{5}$. As expected, $\mathcal{D}_{\xi}^{3}$ is lower bounded by the  fixed mode ambiguity sets  $\mathcal{D}_{\xi}^{2}$ and $\mathcal{D}_{\xi}^{4}$, and upper bounded by $\mathcal{D}_{\xi}^{5}$.

\subsubsection{Reliability}
Using the solutions we generated, we run out-of-sample test with 20 samples of 5000 wind forecast errors to evaluate the joint reliability of each optimal solution. We define the joint reliability as the percentage of wind forecast errors for which all chance constraints are satisfied. Then, we compare the reliability results with our pre-defined probability level ($1-\epsilon=95\%$). The results are summarized in Table~\ref{tab: relia}.
\begin{table*}
\centering
\caption{Joint reliability ($\%$) for $1-\epsilon=95\%$}
\label{tab: relia}
\begin{tabular}{l|c|cccccc|cccc|c|c}
\hline
\hline
       \multirow{2}{*}{Full pool \hspace{0.17cm}}           & \multirow{2}{*}{$\mathcal{D}_{\xi}^{1}$} & \multicolumn{6}{c|}{$\mathcal{D}_{\xi}^{2}$} & \multicolumn{4}{c|}{$\mathcal{D}_{\xi}^{3}$}&  \multirow{2}{*}{$\mathcal{D}_{\xi}^{4}$} & \multirow{2}{*}{$\mathcal{D}_{\xi}^{5}$} \\
       &&  M$1$ & M$2$ & M$3$ & M$4$& M$5$& M$6$\ & $\Xi_1$ & $\Xi_2$ & $\Xi_3$ & $\Xi_4$ & &  \\ \hline\hline
 min &   99.78& 98.22  & 98.08 & 97.72  & 98.14  & 98.02  & 98.06 & 98.54  & 98.48 & 98.58  & 98.52& 98.44  & 98.58  \\
                       avg &  99.87&  98.61 &  98.53 &  98.21 & 98.48  & 98.47  &  98.54 &  98.94 &  98.86 &  98.96 & 98.91&  98.81 &  98.96 \\
                       max &   99.94 &  98.84 &  98.84 &  98.42 &  98.74 &  98.84 &  98.86 &  99.14 &  99.10 &  99.14 &  99.12 &  99.04 &  99.14 \\ \hline
\end{tabular}

\vspace{0.3cm}
\begin{tabular}{l|c|cccccc|cccc|c|c}
\hline
\hline
       \multirow{2}{*}{Partial pool}           & \multirow{2}{*}{$\mathcal{D}_{\xi}^{1}$} & \multicolumn{6}{c|}{$\mathcal{D}_{\xi}^{2}$} & \multicolumn{4}{c|}{$\mathcal{D}_{\xi}^{3}$}&  \multirow{2}{*}{$\mathcal{D}_{\xi}^{4}$} & \multirow{2}{*}{$\mathcal{D}_{\xi}^{5}$} \\
       &&  M$1$ & M$2$ & M$3$ & M$4$& M$5$& M$6$\ & $\Xi_5$ & $\Xi_6$ & $\Xi_7$ & $\Xi_8$ & &  \\ \hline\hline
 min &   99.46 & 92.32  & 92.16 & 82.60  & 88.20  & 91.58  & 92.42& 95.56  & 95.70 & 95.64  & 95.70& 92.72  & 95.78  \\
                       avg &  99.64 &  93.13 &  93.01 &  83.42 & 88.92  & 92.20  &  93.16&  96.24 &  96.29 &  96.24 & 96.29&  93.48 &  96.40 \\
                       max &   99.78 &  93.68 &  93.62 &  84.20 &  89.54 &  92.80 &  93.58 &  96.64 &  96.64 &  96.62 &  96.64 &  93.80 &  96.78 \\ \hline
\end{tabular}
\end{table*}

We observe that reliability ranking almost always matches the cost ranking. Ambiguity sets $\mathcal{D}_{\xi}^{1}$ and $\mathcal{D}_{\xi}^{5}$ have the most conservative solutions and hence  higher reliability and costs. The reliability of $\mathcal{D}_{\xi}^{3}$ is lower bounded by the reliability of $\mathcal{D}_{\xi}^{2}$ and $\mathcal{D}_{\xi}^{4}$, and upper bounded by the reliability of $\mathcal{D}_{\xi}^{5}$. It also shows robustness against the selection of $N_{data}$ and $N_{bin}$. For the full pool, all ambiguity sets achieve constraint satisfaction above $95\%$. For the partial pool, $\mathcal{D}_{\xi}^{2}$ and $\mathcal{D}_{\xi}^{4}$ fail to meet the threshold, while ambiguity sets with misspecified modes $\mathcal{D}_{\xi}^{3}$, arbitrary modes $\mathcal{D}_{\xi}^{5}$, or no unimodality assumptions $\mathcal{D}_{\xi}^{1}$ achieve constraint satisfaction above $95\%$.

In this example, $\mathcal{D}_{\xi}^{5}$ can be use to approximate $\mathcal{D}_{\xi}^{3}$ since they have similar reliability. However, $\mathcal{D}_{\xi}^{3}$ is less conservative than $\mathcal{D}_{\xi}^{5}$ if $\Xi$ does not include the global worst case mode. Set $\mathcal{D}_{\xi}^{3}$ is also more applicable to multivariate unimodality as $\mathcal{D}_{\xi}^{5}$ is only defined for $\alpha=1$.

\subsubsection{Computational Effort}
Table~\ref{tab: time23} shows the iteration count and computational time for $\mathcal{D}_{\xi}^{2}$ and $\mathcal{D}_{\xi}^{3}$. The problems can be solved within 10 iterations and the computational time grows linearly with the number of iterations. Set $\mathcal{D}_{\xi}^{3}$ requires more iterations than $\mathcal{D}_{\xi}^{2}$. Problems using ambiguity sets $\mathcal{D}_{\xi}^{1}$, $\mathcal{D}_{\xi}^{4}$, and $\mathcal{D}_{\xi}^{5}$ can each be solved in a single run, and each takes less than one second.

\begin{table*}
\centering
\caption{Iteration Count and Computational Time for $\mathcal{D}_{\xi}^{2}$ and $\mathcal{D}_{\xi}^{3}$}
\label{tab: time23}
\begin{tabular}{c|cccccc|ccccc}
\hline
\multirow{2}{*}{Full pool} & \multicolumn{6}{c|}{$\mathcal{D}_{\xi}^{2}$} & \multicolumn{4}{c}{$\mathcal{D}_{\xi}^{3}$}\\
& M1 & M2 & M3 & M4 & M5 & M6 & $\Xi_1$ & $\Xi_2$ & $\Xi_3$ & $\Xi_4$ \\ \hline\hline
Iterations & 4  &  4 &  8 &  8 & 4  &  4 &  9 & 8  &  9 &  6 \\
Time (s) & 16.73  & 16.65  &  40.25 & 39.34  & 17.39  &  16.93 & 33.53  & 31.64  & 33.58  & 19.53 \\ \hline
\end{tabular}

\vspace{0.3cm}

\begin{tabular}{c|cccccc|ccccc}
\hline
\multirow{2}{*}{Partial pool} & \multicolumn{6}{c|}{$\mathcal{D}_{\xi}^{2}$} & \multicolumn{4}{c}{$\mathcal{D}_{\xi}^{3}$}\\
& M1 & M2 & M3 & M4 & M5 & M6 & $\Xi_5$ & $\Xi_6$ & $\Xi_7$ & $\Xi_8$ \\ \hline\hline
Iterations & 4  &  4 &  6 &  7 & 4  &  4 & 9 & 9  &  9 &  9 \\
Time (s) & 16.78  & 17.05  &  27.30 & 33.35  & 16.79  &  16.91 & 34.08  & 36.72  & 36.21  & 36.08 \\ \hline
\end{tabular}
\end{table*}

\section{Conclusion}\label{sec:summary}

In this paper, we proposed a distributionally robust chance constrained optimal power flow formulation considering uncertainty distributions with known moments and generalized unimodality with misspecified modes. We derived an efficient solving algorithm using the separation approach. In each iteration of the algorithm, the problem contains only SOC constraints and hence can be solved with commercial solvers. Using wind forecast errors, we found that the distribution of mode estimates are highly dependent on the data pool size, the data size of each sample, and the number of bins used in the histogram. We tested our approach on a modified IEEE 30-bus system and compared our results to those generated with other ambiguity sets. Without the assumption of unimodality, we obtain overly conservative results as unrealistic distributions are included in the ambiguity set. Considering unimodality, but with fixed mode, the results are highly dependent on the quality of the mode estimate. Considering unimodality with misspecified mode, the results are relatively consistent across different mode supports and the performance is bounded by that of the fixed-mode model and that of the arbitrary-mode model. With univariate unimodality and large mode deviations, the misspecified-mode model can be well approximated by the arbitrary-mode model.

Future work will extend the current results by considering more accurate descriptions of the mode support. For example, we could represent the mode support as a union of disjoint sets that matches the mode profile. Other directions include evaluating the approach on a more realistic system and studying how the current approach works in the cases with other misspecified information such as moments.

\appendices
\section{Convexity and Concavity of \eqref{eq: cc-mode}}\label{append1}
Here we prove the left side of \eqref{eq: cc-mode} is neither jointly convex nor concave in $h$ and $\tau$ through counter examples. We first pick $\alpha=\tilde{R}=1$, $\epsilon=0.05$, and $\tilde{c}=0$ without loss of generality. Then we select two groups of points and calculate the left-side values $v$. Group 1: $[h,\tau,v]=(0.1,2,2.985)$ and $(0.3,3,3.05)$, then the midpoint $(0.2,2.5,3.15)$ has a value higher than line segment value $3.0175$ (concave). Group 2: $[h,\tau,v]=(0.4,11,0.1990)$ and $(0.6,10,-1.5015)$, then the midpoint $(0.5,10.5,-0.6693)$ has a value lower than line segment value $-0.65125$ (convex).

\section{Proof of Lemma~\ref{lem: DRCC1}}\label{append2}
We first check if $\underline{h}\geq0$. If so, we know $\tau_1\in\left[\tau_0,\frac{\alpha+1}{\alpha}\right]$ as $\hat{h}\left(\frac{\alpha+1}{\alpha}\right)=0$ and $\hat{h}$ is decreasing. Hence, we can conduct the golden section search on $[\tau_0,\frac{\alpha+1}{\alpha}]$.

Next, if $\underline{h}<0$, we know $\tau_1>\frac{\alpha+1}{\alpha}$ and we have
$$\hat{h}(\tau_1)<h_2(\tau_1)=\frac{f(\tau_1)}{\sqrt{{\frac{1-\epsilon}{\epsilon}+{f(\tau_1)}^2}}}\tilde{R}.$$
If we further force $h_2(\tau_1)=\underline{h}$, we have $\hat{h}(\tau_1)<h_2(\tau_1)=\underline{h}$ and $\overline{\tau}\in[\frac{\alpha+1}{\alpha},\tau_1]$. The equality $h_2(\tau_1)=\underline{h}$ will always have a solution on $\left[\frac{\alpha+1}{\alpha},\infty\right]$ as $h_2\left(\frac{\alpha+1}{\alpha}\right)=0$ and as $\tau\to\infty$, $h_2(\tau)=-\tilde{R}$.

Next we solve the equality and find
\begin{align*}
&f(\tau_1)=\underline{h}\sqrt{\frac{1-\epsilon}{\epsilon(\tilde{R}^2-\underline{h}^2)}}\Rightarrow\\
&\hspace{2cm}\tau_1=-\left(\underline{h}\sqrt{\frac{1-\epsilon}{\epsilon(\tilde{R}^2-\underline{h}^2)}}-(\alpha+1)\right)/\alpha.
\end{align*}
\section{Proof of Lemma~\ref{lem: DRCC2}}\label{append3}
We have the following relationship because $\tau\geq\tau_0>1$.
$$g^{\prime}(\tau)=\frac{\frac{\alpha}{\epsilon}\tau^{-\alpha-1}}{2g(\tau)}\leq g_2(\tau)=\frac{\frac{\alpha}{\epsilon}\tau^{-\alpha}}{2g(\tau)}.$$
Then, we have the following relationship where $\tau_2$ is the effective upper bound.
$$F_3^{\prime}(\tau_2)\leq F_4(\tau_2)=C_3g_2(\tau_2)-(\tilde{c}+\alpha\underline{h})=0.$$
The last equality will always have solution on $[\overline{\tau},\infty]$ since $F_4(\overline{\tau})\geq F_3^{\prime}(\overline{\tau})>0$ and as $\tau\to\infty$, $F_4(\tau)<0$. By solving the equality, we obtain
$$C_2(\tau_2^{-\alpha})^2+\tau_2^{-\alpha}-(1-\epsilon)=0,$$
where $C_2=\frac{\alpha^2C_3^2}{4\epsilon (\tilde{c}+\alpha\underline{h})^2}$. This is a quadratic equation of $\tau_2^{-\alpha}$ and we find
$$\tau_2=\left[\frac{-1+\sqrt{1+4(1-\epsilon)C_2}}{2C_2}\right]^{-\frac{1}{\alpha}}.$$

\bibliography{ReferencesJournal}
\bibliographystyle{IEEEtran}
\end{document}